\documentstyle[psfig,11pt,amsfonts]{amsart}
\title{Chisini's conjecture for curves with singularities of type $x^n=y^m$}
\author{Sandro Manfredini, Roberto Pignatelli}
\thanks{This research is partially supported by The Excellency Center "Group
theoretic methods in the study of algebraic varieties" of the
National Science foundation of Israel, The Emmy Noether Institute for
mathematice, the Minerva foundation of Germany,  the DFG
Forschungsschwerpunkt ``Globale Methoden in der komplexen Geometrie''
and the EU (under EAGER network).}
\begin{document} 
%
\def\hol{{\mathcal O}}
\def\pimaiusc{{\Pi}}
\def\PP{{\mathbb P}}
\def\ci{{\mathbb C}}
\def\zeta{{\mathbb Z}}
\def\eNNe{{\mathbb N}}
\def\erre{{\mathbb R}}
\def\eep{\varepsilon}
\font\rom=cmsy10
\font\tengt=eufm10 \def\bgo{\hbox{\tengt B}}
\def\sgo{\hbox{\rom S}}
\def\eep{\hbox{$\varepsilon$}}
\font\mat=cmr8
\def\ugu{\hbox{\mat \char61}\hskip1pt}
\font\mate=cmmi8
\def\min{\hbox{\mate \char60}\hskip1pt}
\def\mag{\hbox{\mate \char62}}
\font\mpic=cmsy8
\def\magu{\hbox{\mpic \char21}\hskip1pt}
\def\minu{\hbox{\mpic \char20}\hskip1pt}
\def\piu{\raise 1pt\hbox{\mat \char43}\hskip1pt}
\font\rmp=cmr12
\def\men{\hbox{\rmp \char123}\hskip1pt}
\def\diverso{\ugu\hskip -7pt / \kern 2pt}
\font\matem=cmmi12
\def\clopar{\hbox{\matem \char62}}
\def\opar{\hbox{\matem \char60}}
\def\fre{{\longmapsto}}
\def\vir{,\ldots ,}
\def\implica{\hbox{$\,\Rightarrow\,$}}
\def\es{\exists\,}
\def\perogni{\hbox{$\forall\:$}}
\def\IM{\rm Im}
\def\mapright#1{\smash{\mathop{\longrightarrow}\limits^{#1}}}
\def\mapdown#1{\hbox{\big\downarrow{$\vcenter{\hbox{$\scriptstyle#1$}}$}}}
\let\meno\setminus
\def\ristr#1{\lower4pt\hbox{$|_{#1}$}}
\font\cme=line10
\def\modulo#1#2{\setbox1=\hbox{$#2$}\setbox2=\hbox{\cme \char44}
 \raise 3pt \hbox{\hbox{$#1$}
 \hskip-.27truecm\lower.17truecm\copy2\hskip-.2truecm\lower.3truecm\copy1}}
\def\unoenne#1#2{\hbox{$#1_1\vir#1_{#2}$}}
\def\cunoenne#1#2{\hbox{$#1_1\cdots#1_{#2}$}}
\makeatletter
\dimen0=1pt
\def\hidehrule#1#2{\kern-#1\hrule height#1 depth#2 \kern-#2 }
\def\hidevrule#1#2{\kern-#1{\dimen0=#1
    \advance\dimen0 by#2\vrule width\dimen0}\kern-#2 }
\def\makeblankbox#1#2{\hbox{\lower\dp0\vbox{\hidehrule{#1}{#2}%
    \kern-#1
    \hbox to\wd0{\hidevrule{#1}{#2}%
      \raise\ht0\vbox to #1{}
      \lower\dp0\vtop to #1{}
      \hfil\hidevrule{#2}{#1}}%
    \kern-#1\hidehrule{#2}{#1}}}}
\setbox0=\hbox to8\dimen0{\vbox to 8\dimen0{}\hss}
\setbox1\hbox{\makeblankbox{0pt}{.2pt}}
\newbox\cvd@box
\setbox\cvd@box\hbox to 11\dimen0{\vbox to
11\dimen0{\box1\vss\hbox{\kern2\dimen0\vrule height3\dimen0
depth0pt width6\dimen0}}\vrule height 8\dimen0 width3\dimen0 depth0pt}
\def\cvd@{\lower3pt\copy\cvd@box}
\def\CVD{\enspace\hbox{\cvd@}}
\makeatother
\newtheorem{teo}{Theorem}[section]
\newtheorem{lemma}[teo]{Lemma}
\newtheorem{rem}[teo]{Remark}
\newtheorem{corollario}[teo]{Corollary}
\newtheorem{PROP}[teo]{Proposition}
\newtheorem{congettura}[teo]{Conjecture}
\newenvironment{DIM}{{\bf Proof.}}{\hfill
                       \CVD\bigbreak}
\newenvironment{DEF}{\refstepcounter{teo}{\bf Definition \theteo}}{\\}
\newenvironment{DEF*}{\refstepcounter{teo}{\bf Definition
                      \theteo }}{\relax}
\newenvironment{equaz}{$$\refstepcounter{teo}}{\eqno{\rm
(\theteo)}$$
}
\begin{abstract}
This paper is devoted to a very classical problem that can be
summarized as follows: let $S$ be a non singular compact complex surface, 
$\pi: S \rightarrow \PP^2$ a finite morphism having simple branching, 
$B$ the branch curve: 
then (cf. \cite{fulton2}) ``to what extent does $B$ determine 
$\pi: S \rightarrow \PP^2$''?\\
The problem was first studied by Chisini 
(\cite{chis}) who proved that $B$ determines $S$ and $\pi$, 
assuming $B$ to have only nodes and cusps as singularities, the degree
$d$ of $\pi$ to be greater than $5$, and a very strong hypothesis on the
possible degenerations of $B$, and posed the question if the first or
the third hypothesis could be weakened.\\
Recently Kulikov (\cite{kul}) and Nemirovski (\cite{nemi}) proved the
result for $d \geq 12$, and $B$ having only nodes and cusps as singularities.
\\
In this paper we weaken the hypothesis about the singularities of
$B$:
we generalize the theorem of Kulikov and Nemirovski 
for $B$ having only singularities of type $\{x^n=y^m\}$,
in the additional hypothesis of smoothness for the ramification divisor
(automatic in the ``nodes and cusps'' case).
Moreover we exhibit a family of counterexamples
showing that our additional hypothesis is necessary.
\end{abstract}
\maketitle
\section{Introduction.}\label{intro}
As introduced in the abstract, in this paper we study finite morphisms having
simple branching over curves with singularities of type
$\{x^n=y^m\}$.\\
In order to state the problem and our results, we
need to introduce a little bit of notations.\\
\begin{DEF}\label{NGC}
A {\it normal generic cover} is a finite holomorphic map 
$\pi:S\longrightarrow\ci^2,$ which is an analytic
cover branched over a curve $B$ such that $S$ is a connected normal surface and
the fiber over a smooth point of $B$ is supported on deg$\pi\men1$
distinct points.
\end{DEF}
Two normal generic covers $(S_1,\pi_1),$
$(S_2,\pi_2)$ with the same branch locus $B$ are called (analytically)
{\it equivalent} if
there exists an isomorphism $\phi:S_1\rightarrow S_2$ such that
$\pi_1=\pi_2\circ\phi.$ \\
The main interest for generic covers comes from the well known fact
that, by the Weierstrass preparation theorem, given an analytic surface $S
\subset \ci^n$, a generic projection $S \stackrel{\pi}{\rightarrow}
\ci^2$ is (at least locally, in order to insure deg $\pi < \infty$) a 
normal generic cover branched over a curve (see \cite{GR}). \\
A standard way to study generic covers, is the following: given
a generic cover $\pi: S \rightarrow \ci^2$ with branch curve $B$, 
one defines the {\it monodromy homomorphism} $\rho:\pi_1 (\ci^2 \meno B)
\rightarrow \sgo_{{\rm deg}\;\pi}$ as the action of this fundamental
  group on the fiber of a regular value.\\
The pair $(B, \rho)$ gives the ``building data'' of the cover: one can
reconstruct the cover from $(B,\rho)$ (cf. \cite{Grauert}).\\
Despite of the explicit construction, to understand the singularity of
the cover from the building data is a very difficult problem
(except in specific cases). It is, for example, still an open problem
to classify all the possible ``building data'' coming from smooth surfaces.\\
In \cite{mapi} we give a complete classification of the normal generic
covers bran\-ched over irreducible curves of type $\{x^n=y^m\}$ in terms
of what we called there ``monodromy graphs'': we will recall briefly
the definition of monodromy graphs and the above result
 in section 1. Let us point out that, according
to the Puiseux classification, this class of singularities is a
natural first step for a complete classification.\\
Our first result (to which is devoted section 1), is a ``more
friendly'' classification theorem, that will be crucial in the
following sections.\\
Let
$h,k,a,b$ be positive integers with $(h,k)=1$, and consider the surface 
$S_{h,k,a,b}$ in $\ci^4$ defined by the equations
$hz^k+kw^h-(h+k)x^a=zw-y^b=0$. Let $F:S_{h,k,a,b} \rightarrow \ci^2$ be
the projection on the $(x,y)$-plane. 
\begin{teo}\label{classeq0}
The map $F:S_{h,k,a,b} \rightarrow \ci^2$ is a generic cover 
branched over $x^{a(h+k)}=y^{bhk}$ of degree $h+k$.\\
Conversely, up to exchanging $x$ and $y$,
every generic cover $\pi:S\rightarrow\ci^2$ of
degree $d\geq 3$ branched over $\{x^n=y^m\}$, with $(n,m)\ugu1$,
is equivalent to one of the previous maps.
\end{teo}
In section $2$ we consider the ``global'' case of projective generic covers.
\begin{DEF}\label{PGC}
A {\it projective generic cover} is a finite morphism  
$\pi:S\longrightarrow\PP^2,$ branched over an irreducible curve $B$ such that
$S$ is an irreducible projective surface and
the fiber over a smooth point of $B$ has cardinality deg$\pi\men1$.
\end{DEF}
This is the same as requiring that $\pi^*(B)=2R+C,$ with $R$ irreducible
and $C$ reduced, and that $\pi\ristr{R}:R\rightarrow B$ is 1:1 over
smooth points of $B$.\\ 
As in the previous case, for each irreducible projective surface $S$,
a generic projection
$\pi:S\rightarrow\PP^2$ is a projective generic cover branched over a
(projective plane) curve $B$.\\ 
We say that a projective generic cover is {\it smooth}
if the surface $S$ and the
ramification divisor $R$ are non-singular.\\
Actually, when $S$ is non-singular, a ``general'' generic projection
has ramification divisor $R$ non-singular. 
Let us point out that, if $B$ has
only nodes and cusps as singularities, $R$ is automatically smooth.\\
Again, we will consider projective generic covers up to analytic equivalence: 
$(S_1,\pi_1),$
$(S_2,\pi_2)$ with the same branch locus $B$ are {\it equivalent} if
there exists an isomorphism $\phi:S_1\rightarrow S_2$ such that
$\pi_1=\pi_2\circ\phi$.\\ 
Chisini's conjecture asserts the following (see \cite{chis})
\begin{congettura}[Chisini]
Let $B$ be the branch locus of a smooth projective cover $\pi:S\rightarrow\PP^2$ of degree
deg$\pi\magu 5$. Then $\pi$ is unique up to equivalence.
\end{congettura}
In other words, if $S$ is smooth and the degree high enough, the curve $B$
determines the cover.\\
In fact, Chisini proved the result in the above
mentioned additional
hypothesis that the branch curve $B$ has only nodes and cusps as
singularities, and that $B$ has some particular degeneration.
In the same paper, he posed the
question if this two last hypothesis could be weakened.\\
The bound for the degree of $\pi$ is needed according to a
counterexample, due to
Chisini and Catanese (see \cite{cata}) of a sextic curve with 9 cusps which is
the branch curve of 4 non equivalent smooth projective covers, three
of them are of degree 4 and one is of degree 3.\\
Recently, V.S. Kulikov (see \cite{kul}) developed a new approach
proving Chisini's conjecture for curves with only nodes and cusps as
singularities, 
and the additional hypothesis that the degree of $\pi$
is greater than a certain function of the degree, genus and number of
cusps of the branch locus.
After that, S. Nemirovski (see \cite{nemi}), using the
Bogomolov-Miyaoka-Yau inequalities, found a
uniform bound, 12, for Kulikov function.\\ 
Putting the two results together we have the following theorem:
\begin{teo}[\cite{kul}, \cite{nemi}]\hspace*{-5pt}
Let $B$ be the branch locus of a smooth projective cover
$\pi:S\rightarrow\PP^2$ of degree
deg$\pi\magu 12$, with only nodes and cusps as singularities.
Then $\pi$ is unique up to equivalence.
\end{teo}
In section \ref{chisini}, we use theorem \ref{classeq0} in order to improve the
previous results as follows:
\begin{teo}
Let $B$ be the branch locus of a smooth projective generic cover
$\pi:S\rightarrow\PP^2$ 
having only singularities of type $x^{n_i}=y^{m_i}$.
Then, if 
$${\rm deg}\pi>\frac{4(3d+g-1)}{2(3d+g-1)-\sum_{i=1}^r({\rm min}(m_i,n_i)-\gcd{(m_i,n_i)})}$$
where $2d=$deg$B$ and $g=g(B)$ is its genus, then $\pi$ is unique.
\end{teo}
\begin{teo}
In the above hypothesis, if deg$\pi\magu 12$ then $\pi$ is unique.
\end{teo}
Finally, in section \ref{dieci},  we will construct a family of
projective generic covers and we will show that the hypothesis of
smoothness for $R$ is necessary, finding pairs of non equivalent
projective generic covers of arbitrarily large degree
having the same branch curve. More precisely we prove (we defer the 
definitions of $\bar f_i,\bar g_j$ to section \ref{dieci})
\begin{PROP}
Let $t \in \eNNe$, $t\geq 1$, 
$B$ be the projective plane curve given by the equation
$$
\bar g_{4t+1}(x,w)^{2t(2t+1)}=\bar f_{2t(2t+1)}(y,w)^{4t+1}.
$$
Then there are two generic covers $S' \stackrel{\pi'}{\rightarrow}
\PP^2$, $S'' \stackrel{\pi''}{\rightarrow} \PP^2$, with  
$S'$, $S''$ smooth, degrees respectively $4t+1$ and $4t+2$.\\ 
The ramification divisor is singular except in the case $t=1$ and the degree of the cover is 6.
\end{PROP}

{\bf Acknowledgments:}
We would like to thank Prof. Fabrizio Catanese, who was the first to address
us to the subject, for the several useful and interesting conversation
on the topic of the classification of generic covers.\\
We are indebted also with Prof. Victor S. Kulikov who suggested us
the possible applications of our research and
Prof. Mina Teicher who partially supported this research hosting the
first author at Bar Ilan University (Israel).
\section{Equations.}\label{otto}
Consider the following surface $S_{h,k}$ in $\ci^4$ ($S_{h,k,1,1}$ in
  the introduction)
\begin{equaz}\left\{
\begin{array}{l}
hz^k+kw^h=(h+k)x\\
zw=y
\end{array} 
\right.\label{essehk}
\end{equaz}
where $1\minu h\min k$ are coprime integers.\\
The jacobian matrix is
$$\left(\begin{array}{cccc}
h+k & 0 & hkz^{k-1} & hkw^{h-1} \\
0 & 1 & w & z
\end{array}\right)$$
from which we see that $S_{h,k}$ is smooth and we can choose $z,w$ as
local coordinates near $(0,0,0,0)$ for $S_{h,k}$.\\
Consider the map $F_{h,k}:S_{h,k}\longrightarrow\ci^2$ 
which is the restriction to $S_{h,k}$ of the projection of $\ci^4$
on the $(x,y)$-plane.
\begin{PROP}
$F_{h,k}$ is a normal generic cover of degree $h+k$ branched over the
curve $x^{h+k}=y^{hk}$
\end{PROP}
\begin{DIM}
We have that $F_{h,k}^{-1}(0,0)=(0,0,0,0)$ and one can easily check that
the degree of $F_{h,k}$ is $h\piu k$.\\ 
The equations of the ramification divisor $R$ in the local
coordinates $(z,w)$ of $S_{h,k}$ are given by the vanishing of the
determinant of the submatrix of
the jacobian matrix
$$\left (
\begin{array}{cc}
hkz^{k-1} & hkw^{h-1}\\
w & z
\end{array}
\right)$$
that is $z^k=w^h$.\\
Substituting into the equations of $S_{h,k}$ in $\ci^4$, we get that the
locus 
defined by the equation 
$y^{hk}=(z^{hk}w^{hk}=)x^{h+k}$ in the $(x,y)$-plane contains the branch
curve $B$.
But this locus is irreducible since $(h,k)\ugu1$ , so we found the equation of
the branch curve.\\
We are left with the ``genericness'' check. Of course (by irreducibility), it is enough to
check it over a smooth point of $B$, and we take the point $(1,1)$.\\
$F_{h,k}^{-1}(1,1)$ is the set of points of the form $(1,1,z,w)$
described by the equations $\{\frac{hz^k+kw^h}{h+k}=zw=1\}$.\\
Then $z\neq 0$, $w=\frac{1}{z}$ and, 
(multiplying by $z^h$), we have to compute the solutions  of
$$\left({\frac{hz+k}{h+k}}\right)^{h+k}=z^h,$$
i.e. the roots of the polynomial
$P(z)=(hz+k)^{h+k}-(h+k)^{h+k}z^h$.\\
We have to show that $P$ has exactly $h+k-1$ distinct roots; its first and second derivatives are
$$P'(z)=h(h+k)[(hz+k)^{h+k-1}-(h+k)^{h+k-1}z^{h-1}]$$ 
$$P''(z)=h(h+k)[h(h+k-1)(hz+k)^{h+k-2}-(h-1)(h+k)^{h+k-1}z^{h-2}].$$
But $P(z)=P'(z)=0$ implies
$$(hz+k)(h+k)^{h+k-1}z^{h-1}=(h+k)^{h+k}z^h$$
and since $0$ is not a root of $P$, $hz+k=(h+k)z$, i.e. $z\ugu1$.\\
Since $P(1)=P'(1)=0$ but $P''(1)\not=0$, we conclude that $z\ugu 1$ is a
double root of $P$ and all the others are simple roots. 
\end{DIM}
From the proof of previous proposition we get also
\begin{rem}\label{R}
The ramification divisor $R$ is cut (on $S_{h,k}$) by the hypersurface
$z^k=w^h$, while the preimage of the branch locus $B$ is
$2R+C$ where $C$ is the union of the curves cut by the hypersurfaces
$z^k=\alpha w^h$ for $\alpha\diverso 1$ a root of
$P(t)=(ht+k)^{h+k}-(h+k)^{h+k}t^h$.
\end{rem}
Now we introduce the complete class of covers we need for our
classification theorem.

Consider the pullback of $F_{h,k}$ under the base change given by the map
$$f_{a,b}:\ci^2 \rightarrow\ci^2,\ \ f_{a,b}(x,y)=(x^a,y^b).$$
We obtain the surface $S_{h,k,a,b}$ of equations
\begin{equaz}\left\{
\begin{array}{l}
hz^k+kw^h=(h+k)x^a\\
zw=y^b
\end{array} 
\right.\label{essebar}
\end{equaz}
and the map $F_{h,k,a,b}: S_{h,k,a,b} \rightarrow \ci^2$ given by the
two coordinates $(x,y)$.

Now we can introduce the main result of this section.
\begin{teo}\label{classeq}
The maps $F_{h,k,a,b}$ are generic covers of degree $h+k$, 
branched over $x^{a(h+k)}=y^{bhk}$ 

Conversely, up to exchanging $x$ and $y$,
every generic cover $\pi:S\rightarrow\ci^2$ of
degree $d\geq 3$ branched over $\{x^n=y^m\}$, with $(n,m)\ugu1$,
is equivalent to one of the previous maps.
\end{teo}
The first part of the statement is the following lemma.
\begin{lemma}
The maps $F_{h,k,a,b}$ are normal generic covers of degree $(h+k)$ branched
over the curve $x^{a(h+k)}=y^{bhk}$.
\end{lemma}
\begin{DIM}
The statement comes from previous proposition using the base change map 
$f_{a,b}$. The normality of $S_{h,k}$ implies the normality of
$S_{h,k,a,b}$ by theorem 2.2 of \cite{mapi}.
\end{DIM}
In order to prove the second part, we use the well known fact 
(\cite{Grauert}) already mentioned in the introduction that the
pair (branch curve $B$, monodromy homomorphism) determines the
cover. We will introduce now precisely the monodromy homomorphisms and 
the monodromy graphs that
represent them, in term of which we gave in
\cite{mapi} a classification theorem for generic covers branched
over irreducible curves of type $\{x^n=y^m\}$, result that we will finally
briefly recall.\\ 
Let $(S,\pi)$ be a normal generic cover of degree deg$\pi=d$.\\
Every element in the fundamental group $\pi_1(\ci^2\meno B)$ of the set of
regular values of $\pi$, induces a permutation of the $d\ugu$deg$\pi$
points of the fiber over the base point, thus a homomorphism
$\varphi:\pi_1(\ci^2\meno B)\longrightarrow\sgo_d$, called the {\it monodromy} of the cover.
The ``generic'' condition means that for each geometric loop (i.e. a
simple loop around a 
smooth point of the curve) its monodromy is a transposition.
The homomorphisms with this property are called {\it generic} monodromies.\\
So, in order to classify generic covers $\pi:S \rightarrow \ci^2$ of degree
$d$, with $S$ a normal surface, branched over some curve $B$, one need
to classify generic monodromies $\varphi:\pi_1(\ci^2\meno
B)\rightarrow\sgo_d$.\\ 
We did it (for curves $B$ of type $\{x^n=y^m\}$) in \cite{mapi},
representing the monodromy of a normal generic cover of
degree $d$ branched on the curve $\{x^n=y^m\}$ by a labeled graph
$\Gamma$, called {\it monodromy graph}. We will denote by $Gr_{d,n}$
the set of all (isomorphism classes) of graphs with $d$ vertices and
$n$ labeled edges.\\ 
Note that the monodromies of equivalent generic covers differ by an inner
automorphism of $\sgo_d$, so we will say that two monodromies
$\varphi_1,\varphi_2:\pi_1(\ci^2\meno B)\rightarrow\sgo_d$ are {\it
  equivalent} if there exists $\sigma\in\sgo_d$ such that 
$$\varphi_1(\gamma)=\sigma\varphi_2(\gamma)\sigma^{-1}$$
for all $\gamma\in\pi_1(\ci^2\meno B)$.\\
The representation is done as follows: let $\varphi:\pi_1(\ci^2\meno
B)\rightarrow\sgo_d$ be a generic monodromy;
if $\unoenne{\gamma}n$ is a set of geometric loops that ge\-ne\-ra\-tes
$\pi_1(\ci^2\meno B) \cap \{y=1\}$ (in particular they generate
$\pi_1(\ci^2\meno B)$, see, e.g., \cite{oka}, \cite{mapi} for more a
detailed description of this fundamental group), 
we write $d$ vertices labeled $\{1, \ldots, d\}$. Now $\sgo_d$ acts
naturally on the set of our vertices, and then, $\forall i \in
\{1, \ldots, n\}$, we draw the edge labeled $i$ between the two points
exchanged by $\varphi(\gamma_i)$. Finally we have to delete the
labeling of the vertices (this corresponds to consider $\varphi$ up to the equivalence relation
introduced above).\\
Remark that the monodromy graph does not carry all the
informations needed to reconstruct the cover: $\Gamma$ has $n$ edges,
but we lost $m$.\\
For a fixed $\Gamma\in Gr_{d,n}$, we say that $m$ is {\it
  compatible} with $\Gamma$ if $\Gamma$ defines a normal generic cover 
branched over $x^n=y^m$.\\
Then, a pair $(\Gamma\in Gr_{d,n},m)$ with $m$ compatible with $\Gamma$, determines
the cover.\\
Finally, one notes that this construction is not symmetric in the two
variables $x,y$. So, simply exchanging $m$ and $n$, one gets a natural 
involution that sends compatible pairs $(\Gamma\in Gr_{d,n},m)$ in
compatible pairs $(\Gamma'\in Gr_{d,m},n)$: we call this operation
``duality''.\\ 
We need a last definition:\\
\begin{DEF}\label{pol}
 A polygon with $d$ vertices, valence $a$ and increment $j$,
 with $j$ and $d$ relatively prime, is a graph with $n\ugu ad$, 
$d$ vertices, such that $\forall s,t$ the edges labeled $s$ and $t$ have

two vertices in common if and only if $s\men t\ugu\lambda d$

one vertex in common if and only if $s\men t\ugu\lambda d \piu
        j$ or  $s\men t\ugu\lambda d\men j$

no vertices in common otherwise.
\end{DEF}
This complicated definition is in fact probably better explained by
the following example.
\begin{center}\hbox{}
\psfig{figure=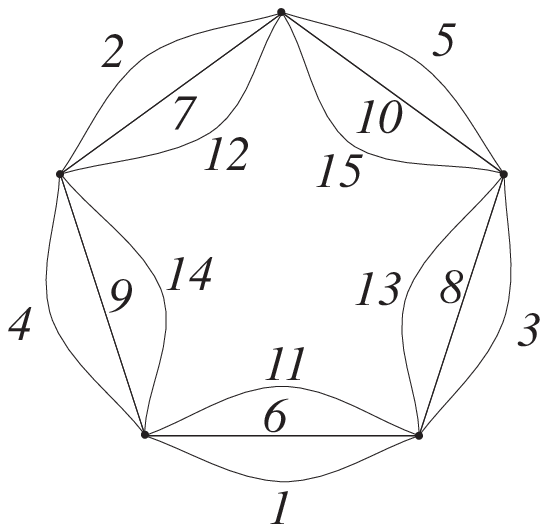}
A polygon with 5 vertices, valence 3 and increment 2.
\end{center}
Now we are able to introduce the main result of \cite{mapi}.
\begin{teo}\label{classgraf}
The monodromy graphs for generic covers $\pi : S \rightarrow \ci^2$ of
degree $d \geq 3$ branched over the curve $\{x^n=y^m\}$, with $(n,m)\ugu1$,
are the following:
\begin{enumerate}
\item ''Polygons'' with $d$ vertices, valence $\frac{n}{d}$ (or
  $\frac{m}{d}$) and increment  $j$, with $(j,d)\ugu1$, $j < \frac{d}{2}$,
  $j(d\men j)|m$ (resp. $j(d\men j)|n$) . Moreover, $d$ must
  divide $n$ (resp. $m$).
\item ''Double stars'' of type $(j,d\men j)$ and valence $\frac{n}{j(d\men j)}$
(or
  $\frac{m}{j(d\men j)}$), with $(j,d)\ugu1$, $j < \frac{d}{2}$, $j(d\men j)|n$
  (resp. $j(d\men j)|m$).  Moreover, $d$ must divide $m$ (resp. $n$).
\end{enumerate}
Duality takes graphs of type $1$ in graphs of type $2$, and vice-versa.
\end{teo} 
We skip here the definition of the double stars (cf. \cite{mapi}), that
  we do not need. \\
Shortly, in theorem \ref{classgraf}, we have shown
that generic covers bran\-ched over an irreducible curve of type $\{x^n=y^m\}$
are classified by pairs (polygon in $Gr_{d,n}$, $m$ multiple of
  $j(d-j)$), up to exchanging $x$ and $y$.\\
Let us recall
that the above pairs describe generic covers also when the hypothesis
$(n,m)= 1$ fails, but in this case we have examples of covers that
can't be described in this way (with monodromy graphs of different type).\\
In view of theorem \ref{classgraf}, in order to prove the remaining
  part of theorem \ref{classeq}, we only need the following
\begin{PROP}\label{monsurf}
The normal generic cover branched over $x^{an}=y^{bm}$ associated to
the polygon with $n$ edges, increment $h$ and valence $a$ is the cover
$F_{h,n-h,a,b}$. 
\end{PROP}
\begin{DIM}
We have to compute the monodromy graphs of the covers $F_{h,k,a,b}$.\\
Let us start by considering the case $a=b=1$, i.e. the covers $F_{h,k}$.\\
$F_{h,k}$ is a normal generic cover branched over $B=\{x^n\ugu y^m\}$ with 
$n\ugu h\piu k$ and $m\ugu hk$. Notice that the assumption $(h,k)=1$ implies
$(n,m)=1$.\\
By theorem \ref{classgraf} the monodromy graph $\Gamma$ is, up to
exchanging $x$ and $y$, a polygon. In fact we do not need to exchange
$x$ and $y$: otherwise we would have $d|m$, while deg$F_{h,k}\ugu n$,
and $(n,m)=1$.\\
So $\Gamma$ has to be a polygon of valence 1 ($d=n$),
and some increment $h'$. 
Set $k'\ugu n\men h'$.\\
By \cite{mapi}, corollary 4.2, the smoothness of $S_{h,k}$ forces $m=h'k'$
(the minimal compatible integer for $\Gamma$).\\
But now $h'+k'=h+k$ and $h'k'=hk$, then $\{h',k'\}=\{h,k\}$.\\
Summing up we proved that the monodromy graph of $F_{h,k}$ is a polygon
with valence $1$ and increment $h$ (or $k$). Of course, the corresponding $m$ 
is $hk$.\\
Now remark that, $\forall{a,b}$, $F_{h,k,a,b}$ can be
obtained by fiber product from $F_{h,k}$  and the map $f_{a,b}:\ci^2
\rightarrow \ci^2$ defined by $f_{a,b}(x,y)=(x^a,y^b)$. As shown in 
\cite{mapi}, this fiber product acts on the ``building data''
of the cover multiplying the valence by $a$, and the compatible $m$
by $b$. So the corresponding mo\-no\-dro\-my graph is a polygon with $d=h+k$
vertices, valence $a$, increment $h$ (or $k$).\\
Conversely, the cover associated to a pair (``polygon with $n$ edges,
valence $a$, and increment $h$'', $m$) is 
$F_{h,n-h,a,\frac{m}{h(n-h)}}$, as stated.
\end{DIM}
This concludes the proof of theorem \ref{classeq}.\\
One immediately gets the following corollary, whose first statement
completes corollary 4.2 in \cite{mapi}. 
\begin{corollario}\label{smoothness}
The cover $F_{h,k,a,b}$ is smooth
$\iff$ $a\ugu b\ugu1$ or $h\ugu b\ugu1$.\\
The cover and the ramification divisor are both smooth
$\iff$ $h\ugu a\ugu b \ugu 1$.
\end{corollario}
\begin{DIM}
The first statement comes by the equations \ref{essebar}, whence
the second can be easily checked in local coordinates as in remark \ref{R}.
\end{DIM}
In the following section we will use the following consequence:
\begin{corollario}\label{solopol}
Let $n$ and $m$ be coprime integers.\\
There exists a non-singular normal generic cover
$\pi:S\longrightarrow\ci^2$ branched over $x^n=y^m$ for which the 
ramification divisor is non-singular if and only if $|m - n|=1$, or
$d=2$, $n=1$.\\
In the first case the cover is unique of degree $d=$max$(m,n)$ and 
its monodromy graph is the polygon with $d$ edges,
increment $1$ and valence $1$.\\
In the second case the cover is given by the projection on the
$x,y$-plane of the surface $z^2=x-y^m$.
\end{corollario}
\begin{DIM}
For $d\geq 3$, by previous corollary, we have only the covers having
as monodromy graph the polygon with $d$ edges, increment $1$ and
valence $1$.\\
For $d=2$, the remark that for every curve $\{f(x,y)=0\}$ there is exactly
one double cover given by projection on the $x,y$-plane of the surface
$z^2=f$, gives immediately the result.
\end{DIM}
We conclude this section with a direct computation of the monodromy graph
associated to $\pi=F_{h,k,a,b}$, although
we don't need it in the rest of the paper: the uninterested reader can skip directly to the next
section.\\
We refer to \cite{mapi} for the notations used in what follows.\\
In order to see how the minimal standard generators act on the preimage of $(1-\eep,1)$, we
examine the inverse image of the path $(\lambda\beta,1)$ for $0\minu\lambda\minu1$
in the $(x,y)$-plane, where $\beta^{h+k}=1$.\\
Since $zw=1,$ we can substitute for $w$ in the first equation \ref{essehk} obtaining
\begin{equaz}
hz^{h+k}-(h+k)z^h\lambda\beta+k=0\label{eqdc}
\end{equaz}
We claim that if $\lambda\diverso 0,1$, then $z^{h+k}$ is real if and only if
$h\piu k$ is odd and $z=s\beta^{-1}$ with $s$ negative.\\
Indeed, if $z^{h+k}$ is real, then $z=s\beta^{-1}$ for some real $s$ and $s$
is a zero of the real polynomial function
$$f(s)=hs^{h+k}-(h+k)\lambda s^h+k.$$
Since $f'(s)=h(h+k)s^{h-1}(s^k-\lambda)$, $f$ will have $s\ugu0$ as critical point
if $h\mag 1$ and one other critical point $s\ugu\sqrt[k]{\lambda}$ if $k$ is odd,
or two other critical points $s\ugu\pm\sqrt[k]{\lambda}$ if $k$ is even.\\ 
If $s^k\ugu\lambda$, then $f(s)=k(1-\lambda s^h)$ which is strictly positive because either
$s\min0$ or $0\min\lambda,s\min 1$. Thus $f$ has only strictly positive critical values,
hence, it has at most one zero $s_0$ and if $f$ does have a zero, then $h+k$ is odd and
$s_0\min 0$.\\
If $\lambda\ugu1$, the same argument shows that $z^{h+k}$ is real if and only if
$s\ugu1$ or $h\piu k$ is odd and $s\min0$.\\
Note that, since $(h,k)=(h,h\piu k)=(k,h\piu k)=1,$ the equations $z^k=\beta,$ $z^k=\beta^{-1}$
for $\beta\diverso1$, $\beta^{h+k}\ugu1$, have a unique common solution:
namely $z\ugu\beta^s$ where $sk\equiv\men sh\equiv 1{\rm\ mod\ }h\piu k$.\\
Now, if $z_0$ is a root of \ref{eqdc} with $\beta\ugu1$, then $z\ugu\frac{z_0}{\beta^s}$
is a root of \ref{eqdc}, so we may restrict to the case $\beta\ugu1$.\\
Note that if $\lambda\ugu0$ then $z^{h+k}=-\frac kh$, while if $\lambda\ugu\beta\ugu1$
\ref{eqdc} has $z\ugu1$ as double root, a real negative root if $h\piu k$ is odd and no
other real roots.\\
Set $\beta_0\ugu\sqrt[h+k]{h+k}e^{-i\frac{\pi}{h+k}}$ and $\alpha_0\ugu e^{i\frac{2\pi}{h+k}}$.
Then, each component of $F_{h,k}^{-1}(\lambda,1)$ will start from one of the points $\alpha_0^r\beta_0$
(each component from a different point) for $r\ugu0\vir h\piu k-1$. Call $c_r$ the component of
$F_{h,k}^{-1}(\lambda,1)$ which starts from $\alpha_0^{r-1}\beta_0$.\\
Then, $c_1$ is contained in the region $-\pi<(h+k)arg(z)<0$, $c_2$ is contained in the region
$0<(h+k)arg(z)<\pi$ and they both have $z\ugu1$ as ending point. Also, $c_{h+k+3-r}=\overline{c_r}$
for $3\minu r\min\frac{h+k+3}2$ are complex conjugated paths and $c_r$ must be contained in one
of the two regions $(2r-4)\pi<(h+k)arg(z)<(2r-3)\pi$ or $(2r-3)\pi<(h+k)arg(z)<(2r-4)\pi$
(see picture). Note that if $h\piu k$ is odd $ c_{\left[\frac{h+k}2\right]+2}$
is contained in the negative real half-line.
\begin{center}
\begin{picture}(140,140)(-70,-70)
\put(0,0){\line(1,-3){20}}
\put(0,0){\line(1,3){20}}
\put(0,0){\line(1,0){80}}
\put(0,0){\line(-5,3){60}}
\put(0,0){\line(-5,-3){60}}
\multiput(0,0)(15,9){5}{\line(5,3){10}}
\multiput(0,0)(15,-9){5}{\line(5,-3){10}}
\multiput(0,0)(-15,0){6}{\line(-1,0){10}}
\multiput(0,0)(-5,-15){4}{\line(-1,-3){4}}
\multiput(0,0)(-5,15){4}{\line(-1,3){4}}
\put(-60,0){\circle*{3}}
\put(50,30){\circle*{3}}
\put(50,-30){\circle*{3}}
\put(-17,51){\circle*{3}}
\put(-17,-51){\circle*{3}}
\put(35,0){\circle{3}}
\put(-30,0){\circle{3}}
\put(4,31){\circle{3}}
\put(4,-31){\circle{3}}
\put(48,26){\vector(-1,-2){11}}
\put(48,-26){\vector(-1,2){11}}
\put(-15,49){\vector(1,-1){16}}
\put(-15,-49){\vector(1,1){16}}
\put(-57,0){\vector(1,0){25}}
\put(25,3){$1$}
\put(44,10){$c_2$}
\put(44,-18){$c_1$}
\put(-9,44){$c_3$}
\put(-9,-48){$c_5$}
\put(-55,5){$c_4$}
\put(-70,-12){$\alpha_0^3\beta_0$}
\put(30,35){$\alpha_0\beta_0$}
\put(38,-39){$\beta_0$}
\put(-45,50){$\alpha_0^2\beta_0$}
\put(-45,-53){$\alpha_0^4\beta_0$}
\end{picture}\\
Configuration of the paths $F_{h,k}^{-1}(\lambda,1)$ in case $h\piu k\ugu5$.
\end{center}
Number the points $\unoenne{z}{h+k}$in $F_{h,k}^{-1}(1-\eep,1)$ by the path $c_r$ they belong to.
It is clear that $z_1$ and $z_2$ are near $z\ugu1$ and that the action of $\gamma_1$ exchanges
$z_1$ and $z_2$.\\
Now, to see which is the action of $\gamma_{h+1}$, follow the motion of the points over the path
$((1\men\eep)(1\men t),1)$ and the path $(t(1\men\eep)\alpha_0^h,1)$ for $0\minu t\minu 1$.
Recall that the paths over $(t(1\men\eep)\alpha_0^h,1)$ are obtained from the paths $c_r$ by
multiplying by $\alpha_0^{-sh}\ugu\alpha_0$, thus, the action of $\gamma_{1+h}$ exchanges
$z_2$ and $z_3$.\\
By the same argument the action of $\gamma_{1+rh}$ will exchange $z_{r+1}$ and $z_{r+2}$, where
indices are taken to be cyclical (mod$h\piu k$), i.e. the monodromy graph associated to $F_{h,k}$ is
the polygon with $h\piu k$ edges increment $h$ and valence 1 (see \ref{pol}).
\section{Chisini's Conjecture}\label{chisini}
In this section we will obtain similar results as those in \cite{kul} and \cite{nemi}
for curves with singularities of type $x^n=y^m$.\\
Let $B\subset\PP^2$ be an irreducible curve with only singularities of
type $\{x^n=y^m\}$. In the whole section, for every such a curve, we write 
$$Sing(B)=\{\unoenne{p}{r}\}$$
where locally, near $p_i$,  $\perogni i\ugu1\vir r$, $B$ is equivalent to 
$x^{s_in_i}=y^{s_im_i}$ 
with 
$(n_i,m_i)\ugu1$, and we set $n_i\min m_i$ (unless $n_i=m_i=1$).
\begin{PROP}
Suppose $B$ is the branch curve of a smooth projective generic cover (cf. \ref{PGC})
$\pi:S\longrightarrow\PP^2$, and let $R$ be the ramification locus of
$\pi$.\\
Then, restricted to the preimage of a small neighborhood of
$p_i$, $\pi$ is given by deg $\pi -n_is_i$ connected components, 
$U_1\vir U_{s_i}$, $V_1\vir V_{deg\pi-(n_i+1)s_i}$
such
that, restricted to one of the $U_j$, $\pi$ gives a
generic cover of degree $n_i+1$ branched over one of the $s_i$
local irreducible components of $B$ (different components for different $j$),
while restricted to each $V_k$, $\pi$ is an isomorphism.\\
Moreover, if $n_i \geq 2$, then $m_i=n_i+1$, and (locally) $\pi$
restricted to $U_j$ is equivalent to the cover $F_{1,n_i,1,1}$ for each $j=1\vir s_i$.
\end{PROP}
\begin{DIM}
Since we assumed $R$ non-singular, it is locally irreducible: then
for each $p\in R$ there exists a neighborhood $U\ni p$ such that
$\pi(R\cap U)$ is irreducible and hence $\pi\ristr U$ is a smooth normal
generic cover branched over an irreducible curve. Since the image of an
irreducible curve is still an irreducible curve, the cover splits locally
as disjoint union of covers each branched over one of the (local) irreducible
components of $B$.\\
In order to prove the first part of the statement, we still have to
compute the degrees of the cover restricted to the ``relevant''
components, that come directly by corollary \ref{solopol}.\\
In case $n_1\geq 2$, by the assumption of smoothness
of the surface and of the ramification divisor $R$, corollary \ref{solopol} forces
$m_i=n_i+1=$the (local) degree of the cover: the local equation for
these covers comes from proposition \ref{monsurf}.
\hspace*{1cm}\end{DIM}
\begin{rem}
By the degrees computed in the previous proposition  we have that 
deg $\pi \magu{\rm max}\{s_i(n_i+1)\}$.
\end{rem}
We introduce some notations: let $\pi:S \rightarrow \PP^2$ be a smooth projective
generic cover, $B$ the branch curve, $B^*$ the dual curve to $B$, $R$ the
ramification locus, $C:=\pi^*(B)-2R$. \\
We set
$E:=\pi^*(\hol_{\PP^2}(1))$ (so that $K_S=-3E+R$), $N:={\rm deg} \pi$, $d:=\frac{{\rm deg}B}{2}$, $\delta:={\rm
  deg}B^*$, $g:=g(B)=g(B^*)=g(R)$. With a standard abuse of notation,
we will not distinguish a divisor from the associated line bundle.\\
In order to prove the main theorem of this section, we follow now the
arguments of Kulikov in our more general case. Although some of the
proof of Kulikov works without correction, we decided, for the
convenience of the reader, to repeat also those proof, with the
exception of proposition \ref{chisiricorda}.\\   
We start with some numerical relations.
\begin{lemma}\label{depari}\  
\begin{enumerate}
\item $d\in\eNNe$;
\item $E^2=N$;
\item $(E,R)=2d$;
\item $\delta=4d+2g-2-\sum_{i=1}^r s_i(n_i-1)$.
\end{enumerate}
\end{lemma}
\begin{DIM}
By Hurwitz formula we have:
$$2-2g(E)=2N-{\rm deg}B$$
thus ${\rm deg}B$ is even, and (1) is proved; (2) and (3) are trivial.\\
Using a generic projection onto a line
$$e(B)=2de(\PP^1)-\delta-\sum_{i=1}^r (s_in_i-1)$$
Thus
$$2-2g=e(R)=e(B)+\sum_{i=1}^r (s_i-1)
=4d-\delta-\sum_{i=1}^r s_i(n_i-1)$$
since $R$ is the normalization of $B$ and is obtained by separating locally
the irreducible components of $B$, and this completes the proof.
\end{DIM}
\begin{lemma}\label{Rquadro}
$$R^2=3d+g-1
$$ 
\end{lemma}
\begin{DIM} By genus formula
$$2g-2=(K_S+R,R)= (-3E+2R,R)=-6d+2R^2.$$
\end{DIM}
Since $\delta\magu 0$, by lemma \ref{Rquadro} and lemma \ref{depari},
part (4), we have 
\begin{corollario}\label{diseg}
$$\sum_{i=1}^r s_i(n_i-1)\leq 2g-2+4d< 2R^2 = 2(3d+g-1)$$
\hfill\CVD
\end{corollario}
By Hodge's Index Theorem ($E^2$ is positive by definition) we have that
$$\left|\begin{array}{cc}
E^2 & (E,R)\\
(E,R) & R^2
\end{array}\right|=N(3d+g-1)-4d^2\leq 0$$
which gives the following
\begin{corollario}
$$N\leq \frac{4d^2}{3d+g-1}$$
\hfill\CVD
\end{corollario}
We can compute the invariants of $S$:
\begin{lemma}\label{invariantidiS}
$$K_S^2=9N-9d+g-1
$$
$$e(S)=3N+\delta-4d=3N+2g-2-\sum_{i=1}^r s_i(n_i-1)$$
$$\chi(\hol_S)=N+\frac{3g-3-9d-\sum_{i=1}^r s_i(n_i-1)}{12}$$
\end{lemma}
\begin{DIM}
Since $K_S=-3E+R$ we have $K_S^2=9N-12d+R^2$.\\
Using a generic pencil of lines in $\PP^2$ and its preimage in $S$ we get
$$e(S)=2e(E)-N+\delta$$ and $e(E)=-(K_S+E,E)=2N-2d$.\\
From Noether's formula $12\chi(\hol_S)=K_S^2+e(S)=12N+3g-3-9d-\sum_{i=1}^r s_i(n_i-1)$
and we have done.
\end{DIM}
Note that $\sum_{i=1}^r s_i(n_i-1)$ must be divisible by 3.\\
Assume that there exist two non-equivalent smooth projective generic covers $(S_1,\pi_1)$ and
$(S_2,\pi_2)$ with the same branch curve $B$.\\
Write $N_i=$deg $\pi_i$ and $\pi_i^*(B)=2R_i+C_i$, for $i\ugu1,2$.\\
Let $X$ be the normalization of the fiber product $S_1\times_{\PP^2}S_2.$
Denote by $g_i:X\rightarrow S_i$, $\pi_{1,2}:X\rightarrow\PP^2$ the
corresponding natural morphisms, as summarized in the following
diagram
$$\begin{array}{ccc}
X & \stackrel{g_1}{\longrightarrow} & S_1\\
{\scriptstyle g_2} \downarrow & \searrow 
\! \! \! \! \! \raise .5em\hbox{$\scriptstyle \pi_{1,2}$}
& \downarrow {\scriptstyle \pi_1} \\
S_2 & \stackrel{\pi_2}{\longrightarrow} & \PP^2
\end{array}$$
We have deg$g_1=N_2$ and deg$g_2=N_1$, so that
deg$\pi_{1,2}=N_1N_2$.\\
The following result is proved in \cite{kul}, proposition $2$,
section $2$. Although Kulikov assumes, at the very beginning, that $B$
has only nodes and cusps as singularities, this proof does
not require this hypothesis.
\begin{PROP}\label{chisiricorda}
If $(S_1,\pi_1)$ and $(S_2,\pi_2)$
are not equivalent, then $X$ is irreducible.
\end{PROP}
Let $Y$ be the set of the points $p\in X$ such that $\pi_{1,2}(p) \subset {\rm Sing} B$, and $\pi_1$
and $\pi_2$ restricted to neighborhoods of $g_1(p)$ and
$g_2(p)$ respectively, are normal generic covers with different branch loci.
\begin{lemma}\label{singxiny}
Sing$X\subset Y$
\end{lemma}
\begin{DIM}
If $g_1(p) \not\in R_1$ or $g_2(p) \not\in R_2$, then $p$ is clearly
smooth.\\
At a point $p$ such that $p_1=g_1(p)\in R_1$ and $p_2=g_2(p)\in R_2$
we can choose small neighborhoods $V_i (p_i)\subset S_i$ and 
$U (\pi_{1,2}(p))\subset\PP^2$
such that $\pi_i(V_i)=U$ and both $\pi_1\ristr{V_1}$ 
and $\pi_2\ristr{V_2}$ are equivalent (up to possibly different base
changes) to one of the following:\\
a) (if $n_i=1$, or $\pi_{1,2}(p)$ is a smooth point of $B$) \\
$f_{2,1}:\ci^2\longrightarrow\ci^2$ defined by $(x,y)\fre(x^2,y)$;\\
b) (if $n_i \geq 2$)\\
the projection on the $(x,y)$-plane of the surface in $\ci^4$
$$\left\{
\begin{array}{l}
n_i w+z^{n_i}=(n_i+1) x\\
zw=y
\end{array}
\right .
$$
Suppose that the branch loci of $\pi_1\ristr{V_1}$ and $\pi_2\ristr{V_2}$ are the same.
We have that,
in the first case, $V_1\times_U V_2$ has equations in $\ci^4$
$$\left\{
\begin{array}{l}
x_1^2=x_2^2\\
y_1=y_2
\end{array}
\right .
$$
and the normalization of $V_1\times_U V_2$ is the disjoint
union of two smooth surfaces, namely $x_1=x_2$, $y_1=y_2$ and $x_1=-x_2$,
$y_1=y_2$ in $\ci^4$.\\
In the second case $\tilde V=V_1\times_U V_2$ is the surface in $\ci^6$
$$\left\{
\begin{array}{l}
n_iw_1+z_1^{n_i}=n_iw_2+z_2^{n_i}=(n_i+1) x\\
z_1w_1=z_2w_2=y
\end{array}
\right .$$
which has two irreducible components, namely $\tilde V_+$
$$\left\{
\begin{array}{l}
w_1=w_2\\
z_1=z_2\\
n_iw_1+z_1^{n_i}=(n_i+1) x\\
z_1w_1=y
\end{array}
\right .$$ 
which is isomorphic to $V_i$ via $g_i$, and $\tilde V_-$
$$\left\{
\begin{array}{l}
n_iw_1=z_2(z_2^{n_i-1}+z_2^{n_i-2}z_1+\cdots+z_2z_1^{n_i-2}+z_1^{n_i-1})\\
n_iw_1+z_1^{n_i}=n_iw_2+z_2^{n_i}=(n_i+1) x\\
z_1w_1=z_2w_2=y
\end{array}
\right .$$
which is expressed by $g_1$ (resp. $g_2$) as a normal generic cover of degree $N_2-1$
(resp. $N_1-1$) branched over $C_1$ (resp.$C_2$).\\
$\tilde V_+$ and $\tilde V_-$ are both smooth and intersect in
$g_1^{-1}(R_1)\cap g_2^{-1}(R_2)$. The normalization will be the disjoint union
of these two smooth components.
\end{DIM}
Now suppose $p\in Y$ and let $V_1$ (resp. $V_2$) be the neighborhood
of $g_1(p)$ (resp. $g_2(p)$) as in the definition of $Y$: the branch
loci of $\pi_1\ristr{V_1}$ and $\pi_2\ristr{V_2}$ are different. 
\begin{PROP}\label{singolarita}
$X$ has only R.D.P. as singularities.\\
More precisely, for every point $P \in Y$, if $p_i=\pi_{1,2}(P)$,
$P$ is a point of $X$ of type $A_{m_i-1}$, and these are all the
singular points of $X$.
\end{PROP}
For instance, if $n_i=m_i=1$ (the case of nodes), we get $A_0$, i.e. a
smooth point.\\
\begin{DIM}
By previous lemma, $\pi_{1,2}(P)\ugu p_i$ for some $i$.\\
If $n_i=1$, we can assume the two branch loci to be
$\{x=0\}$ and $\{x+y^{m_i}=0\}$ and we get
$$\left\{
\begin{array}{l}
x=z_1^2\\
z_2^2=x+y^{m_i}
\end{array} 
\right.$$
i.e. 
$z_2^2=z_1^2+y^{m_i}$ that is clearly a singularity of type
$A_{m_i-1}$ (if $m_i=1 \Rightarrow X$ is smooth at $P$).\\
Finally, if $n_i \neq 1$, then $m_i=n_i+1$, and $V_1\times_U V_2$ is the surface in $\ci^4$
$$\left\{
\begin{array}{l}
z_1^{m_i}-m_ixz_1=-(m_i-1) y\\
z_2^{m_i}-m_i\alpha xz_2=-(m_i-1) y\\
\end{array}
\right .$$
with $\alpha^{s_im_i}\ugu1$ but $\alpha^{m_i}\diverso1$, which
is isomorphic to the surface in $\ci^3$ 
$$
z_1^{m_i}-m_ixz_1=z_2^{m_i}-m_i\alpha xz_2
$$
which has a double point at the origin.\\
The Hessian matrix ($m_i \geq 3$), at the origin is
$$\left(
\begin{array}{ccc}
0&-m_i&-m_i\alpha\\
-m_i&0&0\\
-m_i\alpha&0&0
\end{array}
\right )$$
and has rank 2, hence $X$ has in $p$ a singularity of type $A_k$ for some $k\geq 2$.\\
In order to compute $k$, set
$z:=m_i(z_1-\alpha z_2)$. Then, in the coordinate system $(x,z,z_1)$,
our equation can be 
written as $z_1^{m_i}=z(x+f(z,z_1))$ with $f(0,0)=0$, and, setting 
$\overline{x}=x+f$, we find that near the origin the triple
$(\overline{x},z,z_1)$ is still a coordinate system in terms of which 
$V_1\times_U V_2$ has equation $z_1^{m_i}=\overline{x}z$ that is the standard
expression for the singularity $A_{m_i-1}$.
\end{DIM}
Note that, if $P$ is singular for $X$, $g_1^{-1}(R_1)\cap
g_2^{-1}(R_2) \cap (V_1\times_U V_2)=P$. \\
In general, if $D_1$ and $D_2$ are two divisors in a normal surface, we
define $(D_1,D_2)$ (``the greatest common divisor'') as the greatest
divisor contained in both. By the local equations for the ramification
divisor given in remark \ref{R} we notice that the ``singular'' points
in proposition \ref{singolarita} are isolated points for $g_1^{-1}(R_1)\cap
g_2^{-1}(R_2)$, then
\begin{rem}\label{Rliscio}
If $R=(g_1^{-1}(R_1),
g_2^{-1}(R_2))$, $R$ does not intersect Sing $X$ and, by the local
considerations in the proof of lemma
\ref{singxiny}, is smooth and $g_i\ristr R:R\rightarrow R_i$
is a (unramified) double cover.
\end{rem}
Let $F:\tilde X\rightarrow X$ be the resolution of singularities of
$X$, and let $\tilde g_i=g_i\circ F$, $\tilde\pi_{1,2}=\pi_{1,2}\circ F$.
We define $\tilde{R}:=F^{*}(R)$,
$\tilde C_1:=F^*((g_1^{-1}(R_1), g_2^{-1}(C_2)))$, 
$\tilde C_2:=F^*((g_1^{-1}(C_1), g_2^{-1}(R_2)))$.
\begin{PROP}
\begin{enumerate}
\item $(\tilde R,\tilde C_j)=\sum_{i=1}^r s_i(n_i-1)$
\item $\tilde R^2= 2(3d+g-1)-\sum_{i=1}^r s_i(n_i-1)$
\item $\tilde C_1^2=(N_2-2)(3d+g-1)-\sum_{i=1}^r s_i(n_i-1)$
\item $\tilde C_2^2=(N_1-2)(3d+g-1)-\sum_{i=1}^r s_i(n_i-1)$
\end{enumerate}
\end{PROP}
\begin{DIM}
By remark \ref{Rliscio} $R$ does not intersect the singular points of $X$, 
then we can compute the intersections of $\tilde C_1$ and $\tilde R$
in $X$.
By definition, $\tilde C_1$ and $\tilde R$ intersect only at points
of the preimage of $R_2 \cap C_2$, in particular over some singular
point of $B$.\\
But we already noticed that the only points $p\in R$, s.t.
$\pi_{1,2}(p) \in $ Sing $B$ are the points such that $\pi_1$ near $g_1(p)$
and $\pi_2$ near $g_2(p)$ are branched over the same curve, considered in the
proof of lemma \ref{singxiny}.\\ 
Let $p\in X$ be such a point.
$\pi_{1,2}(p)$ is a singular point of $B$, so that there exists $i$ such
that $\pi_{1,2}(p)=p_i$.\\
In case $n_i=1$, $\pi_1$ and $\pi_2$ are locally double covers, so
$C_j$ does not contain $p$, and $p$ does not give contribution to the
intersection number.\\ 
Otherwise, let $V_1$ (resp. $V_2$) be a small neighborhood of $g_1(p)$
(resp. $g_2(p)$) as in lemma \ref{singxiny}. Then, since $g_1\ristr{R}$ (and
also $g_2\ristr{R}$) is an unramified double cover, there are exactly two
points over $p_i$ contained both in $R$ and in the normalization of the
fiber product of $V_1$ and $V_2$, say $P_{i+}$ and $P_{i-}$.\\                 The two points belong to the two components $\tilde V_+$ and
$\tilde V_-$ respectively (see proof of lemma \ref{singxiny}), but since
$C_j$ does not intersect $\tilde V_+$, we may suppose $p = P_{i-}\in
\tilde V_-$.\\
If we rewrite the equations for $\tilde V_-$ we get
\begin{eqnarray*}
w_1&=&\frac{z_2(z_1^{n_i}-z_2^{n_i})}{n_i(z_1-z_2)}\\
w_2&=&w_1+\frac{z_1^{n_i}-z_2^{n_i}}{n_i}\\
x&=&\frac{n_iw_1+z_1^{n_i}}{n_i+1}\\
y&=&z_1w_1.
\end{eqnarray*}
Remark that all the members of these equations 
are polynomials, and
we can take $z_1,z_2$ as holomorphic coordinates for $\tilde V_-$.\\
Now, $R\cap\tilde V_-$ is cut by $w_1=z_1^{n_i}$, $w_2=z_2^{n_i}$
i.e.
$$\left\{\begin{array}{l}
z_2\frac{z_1^{n_i}-z_2^{n_i}}{z_1-z_2}
=n_iz_1^{n_i}\\
z_1^{n_i}-z_2^{n_i}+z_2\frac{z_1^{n_i}-z_2^{n_i}}{z_1-z_2}\!=\!n_iz_2^{n_i}.
\end{array}\right. $$
This implies $z_1^{n_i}=z_2^{n_i}$, i.e. $z_1=\lambda z_2$, with 
$\lambda^{n_i}=1$. But if $\lambda \neq 1$, then
the left members in our equations vanish, and we get $z_1=z_2=0$
(which is not a curve), whence $z_1=z_2$ clearly solve our equations,
so it is the local equation we were looking for.\\
A branch of $F(\tilde C_2)\cap\tilde V_-$ is given (cf. remark \ref{R})
by the equations $\alpha w_1=z_1^{n_i}$, 
$w_2=z_2^{n_i}$, where
$\alpha =\left(\frac{n_i+\alpha}{n_i+1}\right)^{n_i+1}$,
$\alpha \neq 1$, i.e.
$$\left\{\begin{array}{l}
z_2\frac{z_1^{n_i}-z_2^{n_i}}{z_1-z_2}
=\frac{n_i}{\alpha}z_1^{n_i}\\
z_1^{n_i}-z_2^{n_i}+z_2\frac{z_1^{n_i}-z_2^{n_i}}{z_1-z_2}\!=\!n_iz_2^{n_i}
\end{array}\right.$$
from which $z_1^{n_i}=\frac{\alpha (n_i+1)}{n_i+\alpha}z_2^{n_i}=\left(\frac{n_i+\alpha}{n_i+1}\right)^{n_i}z_2^{n_i}$ so that
$z_1=\lambda\frac{n_i+\alpha}{n_i+1}z_2$ with $\lambda^{n_i}=1$.\\
Moreover, setting $t=\frac{n_i+\alpha}{n_i+1}$ (so that $t^{n_i+1}=(n_i+1) t-n_i$),
$\lambda$ must satisfy (by the first equation)
$$(\lambda t)^{n_i-1}+(\lambda t)^{n_i-2}+\cdots+1=\frac {n_i}t$$ 
Hence
$$t^{n_i}-1=(\lambda t)^{n_i}-1=(\lambda t-1)\frac {n_i}t$$
or
$$t^{n_i+1}=(\lambda n_i+1)t-n_i$$
i.e.
$$(n_i+1) t=(\lambda n_i+1)t.$$
Thus $\lambda=1$ and $F(\tilde C_2)\cap\tilde V_-$ is the union of the $n_i-1$
curves $z_1=\frac{n_i+\alpha}{n_i+1}z_2$.\\
Then every component of $\tilde C_j$ intersects $\tilde R$ transversally, and we
conclude $(\tilde R,\tilde
C_j) =\sum_{i=1}^r s_i(n_i-1)$.\\
Let $E_{\tilde{X}}=F^*\pi_{1,2}^*(\hol_{\PP^2}(1))=F^*g_i^*(E_i)$. It is immediate to verify that
$$E_{\tilde{X}}^2=N_1N_2$$ 
$$(E_{\tilde{X}},\tilde R)=4d$$ 
$$(E_{\tilde{X}},\tilde C_1)=2d(N_1-2)$$ 
$$(E_{\tilde{X}},\tilde C_2)=2d(N_2-2).$$
Since the canonical divisor of $\tilde{X}$ is
$F^*K_X$ ($X$ has only rational double points as singularities), we have
$$K_{\tilde{X}}=-3E_{\tilde{X}}+\tilde R+\tilde C_1+\tilde C_2$$
and
$$(K_{\tilde X}+\tilde R,\tilde R)=e(\tilde R)=4g-4$$
then
$$\tilde R^2=6d+2g-2-\sum_{i=1}^r s_i(n_i-1).$$
Since $F^*g_1^*(R_1)=\tilde R+\tilde C_1$ we have
$$N_2 R_1^2=(\tilde R+\tilde C_1,\tilde R+\tilde C_1)$$
from which
$$\tilde C_1^2=N_2(3d+g-1)-6d-2g+2-\sum_{i=1}^r s_i(n_i-1)=(N_2-2)(3d+g-1)-\sum_{i=1}^r s_i(n_i-1).$$
\end{DIM}
Finally we can prove
\begin{teo}
Let $B$ be the branch locus of a smooth projective generic cover
$\pi:S\rightarrow\PP^2$ 
having $r$ singular points of type $x^{n_is_i}=y^{m_is_i}$ with $n_i \leq
m_i$, $(n_i,m_i)=1$.
Then, if 
$${\rm deg}\pi>\frac{4(3d+g-1)}{2(3d+g-1)-\sum_{i=1}^r s_i(n_i-1)}$$
where $2d=$deg$B$ and $g=g(B)$ is its genus, then $\pi$ is unique.
\end{teo}
In the introduction we wrote the statement in a different notation,
that we found better there.\\
\begin{DIM}
Since by corollary \ref{diseg} $\tilde R^2>0$, by Hodge Index Theorem
$$\left|\begin{array}{cc}
\tilde R^2 & (\tilde C_1,\tilde R)\\
(\tilde C_1,\tilde R) & \tilde C_1^2
\end{array}\right|=2(N_2-2)(3d+g-1)^2-N_2(3d+g-1)\sum_{i=1}^r s_i(n_i-1)\leq 0$$
and the same equation is true replacing $\tilde C_1$ by $\tilde C_2$ and $N_2$ by $N_1$.\\
So, we get
$$N_j\leq\frac{4(3d+g-1)}{2(3d+g-1)-\sum_{i=1}^r s_i(n_i-1)}$$
for $j=1,2$.
\end{DIM}
Following an idea of S. Nemirovski (see \cite{nemi}) we may prove the following
\begin{teo}
In the above hypothesis, if deg $\pi\magu 12$ then $\pi$ is unique.
\end{teo}
\begin{DIM}
If $S$ is not an irrational ruled surface of genus $g \geq 2$, it
satisfies the Bogomolov-Miyaoka-Yau inequality 
$$K_S^2\leq 3e(S).$$
From lemma \ref{invariantidiS}
$$K_S^2=9N-9d+g-1$$
$$e(S)=3N+2g-2-\sum_{i=1}^r s_i(n_i-1)$$
so,
$$\sum_{i=1}^r s_i(n_i-1)\leq 3d+\frac 53(g-1).$$
With this inequality, we can estimate the quantity
$$\frac{4(3d+g-1)}{2(3d+g-1)-\sum_{i=1}^r s_i(n_i-1)}\leq\frac{12d+4(g-1)}{3d+\frac 13(g-1)}
=4+\frac{8(g-1)}{9d+g-1}<12$$
If $S$ is an irrational ruled surface, it satisfies 
$$K_S^2\leq 2e(S)$$
then, the same argument shows that
$$\sum_{i=1}^r s_i(n_i-1)\leq\frac{-3N+9d+3(g-1)}2<\frac 32(3d+g-1)$$
thus we get the stronger estimate
$$\frac{4(3d+g-1)}{2(3d+g-1)-\sum_{i=1}^r s_i(n_i-1)}<8$$
\end{DIM}
As a last remark, note that one can rewrite, with the obvious changes,
all the results in \cite{kul}, theorems 3-12.
\section{A family and a counterexample}\label{dieci}
In this section we will describe an interesting family
of projective generic covers branched over a
curve $\bar B$ with singularities of type $x^n=y^m$ that will produce
a counterexample to Chisini's conjecture if we drop the hypothesis
that the ramification divisor is smooth.\\
Let $\bar B\subset\PP^2$ be a plane curve of equation $\bar g(x,w)=\bar
f(y,w)$ where 
$\bar g$ and $\bar f$ are homogeneous polynomials of degree $d$,
of the form
$$\bar g(x,w)=\prod_{i=1}^r(x-\alpha_iw)^{n_i}$$
$$\bar f(y,w)=\prod_{j=1}^s(y-\beta_jw)^{m_j}$$
with $\unoenne{\alpha}r$ and $\unoenne{\beta}s$ mutually distinct.\\
In a neighborhood $U_{i,j}$ of the point
$P_{i,j}=(\alpha_i,\beta_j,1)$, $\bar B$ is analytically 
e\-qui\-va\-lent to $x^{n_i}=y^{m_j}$.\\
Our (open) assumption is that the singular points of $\bar B$ are contained
in the union of lines $\bar g(x,w)=0$, or, if you prefer, in the set
of the $P_{ij}$'s.\\
By a classical result (see \cite{del,ful}), if $\bar B$ is a nodal
curve then $\pi_1(\PP^2\meno\bar B)$ is abelian; since $S_d$ has no
center if $d\magu3$, then, if $\pi_1(\PP^2\meno\bar B)$ is abelian,
there are no projective generic covers of degree $d\magu3$ whose
branch locus is $\bar B$; thus we will suppose that not all $n_i\minu
2$ and not all $m_j\minu2$.\\
Note that $p=(0,1,0)$ does not belong to $\bar B$, thus, in order to
compute $\pi_1(\PP^2\meno\bar B)$, we can use the projection from $p$
onto the $x$-axis.\\ 
More precisely, $\bar B$ intersects transversally the line at infinity 
$w\ugu 0$ in the $d$ smooth points $(1,\xi,0)$ with $\xi^d\ugu1$; then
the line at infinity is not tangent to $\bar B$. This allows us to compute the
fundamental group of the complement of $\bar B$ by computing the fundamental
group of the complement of the affine curve $B$ in the chart $w\neq
0$, as we will do in proposition \ref{pi1complproj}.\\  
Set $g(x)=\bar g(x,1)$ and $f(y)=\bar f(y,1)$ so that $B=\{g(x)\ugu f(y)\}$.\\
In order to compute the fundamental group of the complement of $B$ we can do, without
lost of generality (by a deformation argument as in \cite{oka}), the
following assumptions:
\begin{enumerate}
\item\label{real} $\forall i,j$ , $\alpha_i, \beta_j \in \erre$; 
\item $\alpha_1\min\alpha_2\min\cdots\min\alpha_r$, and $\beta_1\min\beta_2\min\cdots\min\beta_s$;
\item \label{crit} If $\unoenne{\gamma}{s-1}$ are the roots of $f'$ such that
$f(\gamma_i) \neq 0$,  the 
critical values for $f$, $f_1=f(\gamma_1)\vir f_{s-1}=f(\gamma_{s-1})$
are mutually distinct;
\item\label{arba} For a suitable $\eep_0\mag 0$, $\forall x \in
  (\alpha_1\men\eep_0, \alpha_r\piu\eep_0)$,  $|g(x)|< {\rm min}_i\ |f_i|$
\end{enumerate}
Let us point out (in order to justify assumption \ref{crit}) that the roots of $f'$ are those $\beta_j$ for which 
$m_j\magu2$ (with multiplicity $m_j-1$) and the roots of a polynomial
of degree $s-1$ that has, by assumption \ref{real}, $s-1$ distinct
real roots $\unoenne{\gamma}{s-1}$ such that 
$\beta_i < \gamma_i < \beta_{i+1}$.\\
The critical points of the projection from $p$ onto the $x$-axis are
given by the intersection of $B$ with the union of horizontal lines
$\{f'(y)=0\}$.\\
Then the critical values are (some of) the $\alpha_i$ (corresponding
to points $P_{i,j}$) and 
the $d(s\men1)$ distinct points $\delta_{j,h}$ for $h\ugu1\vir d$ and
$j\ugu1\vir s\men1$ where $g(\delta_{j,h})=f_j$ (smooth
points with vertical tangent). 
By assumption \ref{arba}, no $\delta_{j,h}$ is contained in the
interval $[\alpha_1,\alpha_r]$.\\
Choose $\eep\mag 0$ small enough such that, $\forall j$
(resp. $\forall i$), for every $t$ s.t.
$0<|t-\beta_j|\leq \eep$ (resp. $0<|t-\alpha_i|\leq \eep$), $f^{-1}(t)$
(resp. $g^{-1}(t)$) is given
by $m_j$ (resp. $n_i$) distinct points.
We denote by $b_{j,1}\vir b_{j,m_j}$, (resp. $a_{i,1}\vir a_{i,n_i}$)
the points in $f^{-1}(\eep)$ (resp. $g^{-1}(\eep)$) ordered by their
argument.\\ 
We fix now a free basis for $\Pi=\pi_1(\{y\ugu 0\}\meno
\{\alpha_i, \delta_{j,k}\},a_{1,1})$, in terms of which we will describe the
braid monodromy of the projection.
\\
Let $C_{\eep}(z_0)\subset\ci$ be the circle of center $z_0$ and
radius $\eep$; we define by $c_i$ the closed path supported on the
connected component of $g^{-1}(C_{\eep}(0))$ near $\alpha_i$, with
starting point the unique
real point bigger than $\alpha_i$, with the counterclockwise
orientation; $c_i^+$ the ``subpath'' contained in the positive half
plane (imaginary part bigger than $0$), 
 $c_i^-$ the ``subpath'' in the negative half plane.\\
Let $l_i$ (for $i\ugu1\vir r\men1$)
be the (positively oriented) path contained in the real line connecting $c_i\cap\erre$ and
$c_{i+1}\cap\erre$ but not containing any of the $\alpha_j$'s. 
\begin{center}
\begin{picture}(350,100)(-25,-40)
\put(0,0){\line(1,0){218}}
\put(30,0){\circle{30}}
\put(100,0){\circle{30}}
\put(287,0){\circle{30}}
\put(170,0){\circle{30}}
\put(30,0){\circle*{3}}
\put(100,0){\circle*{3}}
\put(170,0){\circle*{3}}
\put(287,0){\circle*{3}}
\put(220,0){$\ldots$}
\put(238,0){\line(1,0){80}}
\put(41,11){\circle*{3}}
\put(41,-11){\circle*{3}}
\put(19,11){\circle*{3}}
\put(19,-11){\circle*{3}}
\put(116,0){\circle*{3}}
\put(154,0){\circle*{3}}
\put(179,13){\circle*{3}}
\put(287,-16){\circle*{3}}
\put(287,16){\circle*{3}}
\put(303,0){\circle*{3}}
\put(271,0){\circle*{3}}
\put(179,-13){\circle*{3}}
\put(26,3){$\alpha_1$}
\put(96,3){$\alpha_2$}
\put(166,3){$\alpha_3$}
\put(283,3){$\alpha_r$}
\put(25,-13){$c_1$}
\put(95,-13){$c_2$}
\put(165,-13){$c_3$}
\put(282,-13){$c_r$}
\put(282,-24){$a_{r,4}$}
\put(282,20){$a_{r,2}$}
\put(60,3){$l_1$}
\put(130,3){$l_2$}
\put(200,3){$l_3$}
\put(245,3){$l_{r-1}$}
\put(39,14){$a_{1,1}$}
\put(39,-19){$a_{1,4}$}
\put(-1,14){$a_{1,2}$}
\put(-1,-19){$a_{1,3}$}
\put(178,16){$a_{3,1}$}
\put(178,-19){$a_{3,3}$}
\put(115,-8){$a_{2,1}$}
\put(302,-8){$a_{r,1}$}
\put(136,-8){$a_{3,2}$}
\put(253,-8){$a_{r,3}$}
\end{picture}
$n_1\ugu n_r\ugu4$, $n_2\ugu1$, $n_3\ugu3$
\end{center}
Let $\omega$ be the small path supported on $c_1$ connecting
$a_{1,1}$ with the base point of $c_1$ (in the clockwise direction).\\
Consider the paths $\unoenne{\rho}r$ based at $a_{1,1}$ defined by
$$
\rho_i=(\omega l_1(c_2^+)^{-1}\cdots l_{i-2}(c_{i-1}^+)^{-1}l_{i-1})c_{i}
(\omega l_1(c_2^+)^{-1}\cdots l_{i-2}(c_{i-1}^+)^{-1}l_{i-1})^{-1}
$$  
($\rho_1=\omega c_1\omega^{-1}$).
\begin{center}
\begin{picture}(350,100)(-25,-40)
\put(0,0){\line(1,0){218}}
\put(30,0){\circle{30}}
\put(100,0){\oval(31,31)[b]}
\put(165,-23){$\rho_3$}
\put(167,14){$\min$}
\put(30,0){\circle*{3}}
\put(100,0){\circle*{3}}
\put(170,0){\circle*{3}}
\put(287,0){\circle*{3}}
\put(220,0){$\ldots$}
\put(238,0){\line(1,0){80}}
\put(41,11){\circle*{3}}
{\thicklines
\qbezier(41,11)(46,7)(46,0)
\put(46,0){\line(1,0){39}}
\put(100,0){\oval(31,31)[t]}
\put(115,0){\line(1,0){39}}
\put(170,0){\circle{30}}}
\put(26,3){$\alpha_1$}
\put(96,3){$\alpha_2$}
\put(166,3){$\alpha_3$}
\put(283,3){$\alpha_r$}
\put(39,14){$a_{1,1}$}
\end{picture}
\end{center}
These are path around the $\alpha_i$'s. To complete the free basis of
$\Pi$, we need some paths around the $\delta_{j,k}$'s.\\
Consider the (real) critical values for $f$, $f_i$, defined before.\\
Let $\omega_i$ be a loop around $f_i$
based at $\eep$ contained in the union of the paths $C_{\eep}(f_i)$
and the real line constructed by the following algorithm: follow the real
line in direction of $f_i$ until you meet the first $C_{\eep}(f_j)$; if $j
\neq i$, follow $C_{\eep}(f_j)$ clockwise until you meet again the real
line, then follow the real line again till a new $C_{\eep}(f_j)$ and
repeat the algorithm; if $i=j$ follow counterclockwise the whole 
$C_{\eep}(f_i)$ and come back to $\eep$ from the way you arrived
(and end the algorithm).
Here you find two examples, were we defined $C_{\eep}^+$ and 
$C_{\eep}^-$ in the natural way as we did for the $c_i$.
\begin{center}
\begin{picture}(350,100)(0,-40)
\put(0,0){\line(1,0){300}}
\put(90,0){\oval(31,31)[t]}
\put(230,0){\oval(30,30)[b]}
\put(160,0){\oval(30,30)[t]}
\put(20,0){\circle*{3}}
\put(90,0){\circle*{3}}
\put(160,0){\circle*{3}}
\put(230,0){\circle*{3}}
\put(300,0){\circle*{3}}
\put(238,0){\line(1,0){80}}
\put(175,0){\circle*{3}}
{\thicklines
\put(230,0){\oval(30,30)[t]}
\put(20,0){\circle{30}}
\put(160,0){\oval(30,30)[b]}
\put(176,0){\line(1,0){39}}
\put(245,0){\line(1,0){39}}
\put(36,0){\line(1,0){39}}
\put(90,0){\oval(31,31)[b]}
\put(105,0){\line(1,0){40}}
\put(300,0){\circle{30}}}
\put(16,4){$f_3$}
\put(86,4){$f_2$}
\put(155,3){$0$}
\put(226,4){$f_4$}
\put(296,4){$f_1$}
\put(175,2){$\eep$}
\put(297,14){$\min$}
\put(17,14){$\min$}
\put(15,-23){$\omega_3$}
\put(295,-23){$\omega_1$}
\end{picture}
$\omega_1=TC_{\eep}(f_1)T^{-1}$
$T=[\eep,f_4-\eep]C_{\eep}^+(f_4)^{-1}[f_4+\eep,f_1-\eep]$\\
$\omega_3=T'C_{\eep}(f_3){T'}^{-1}$
$T'=C_{\eep}^-(0)^{-1}[-\eep,f_2+\eep]C_{\eep}^-(f_2)^{-1}[f_2-\eep,f_3+\eep]$
\end{center}
For every fixed pair $i,h$, we can uniquely lift $\omega_j$ to a
(closed) path $\tilde{\Delta}_{j;i,h}$, based at $a_{i,h}$,
s.t. $g(\tilde{\Delta}_{j;i,h})=f(\omega_j)$, that is in fact a loop
around some $\delta_{j,\bar h}$.\\ 
Finally we define $\Delta_{j;i,h} \in \Pi$ the path based in $a_{1,1}$
obtained conjugating $\tilde{\Delta}_{j;i,h}$ by a path connecting
$a_{1,1}$ and $a_{i,h}$, obtained following the
orientation of each real interval and the reverse orientation of each
circle.\\
The paths $\rho_i$'s, $\Delta_{j;i,h}$'s give clearly a free basis for $\Pi$.\\
Now we can compute $\pi_1(\ci^2\meno B)$ (and  $\pi_1(\PP^2\meno\bar B)$).\\
We can take as generators of $\pi_1(\ci^2\meno B)$ (and of $\pi_1(\PP^2\meno\bar B)$)
a geometric basis $\mu_{j,k}$ (for $j\ugu1\vir s$, $k\ugu1\vir m_j$) of
$\pi_1(\{x=a_{1,1}\}\meno B)\cong F_{d}$ in such a
way that $\mu_{j,1}\vir\mu_{j,m_j}$ are (conjugated to) the ``standard
generators'' of $\pi_1(U_{1,j}\meno B)$ (cf. \cite{mapi}) as in the
following picture.
\begin{center}
\begin{picture}(190,120)(94,-20)
\put(167,2){\makebox(1,2){$(a_{1,1},b_{1,1})$}}
\put(90,0){\makebox(1,2){$(a_{1,1},\beta_1)$}}
\put(119,53){\makebox(1,2){$(a_{1,1},b_{1,2})$}}
\put(59,57){\makebox(1,2){$(a_{1,1},b_{1,i})$}}
\put(156,66){\makebox(6,10){$\mu_{1,2}$}}
\put(164,31){\makebox(6,10){$\mu_{1,1}$}}
\put(26,70){\makebox(6,10){$\mu_{1,i}$}}
\put(90,68){\makebox(6,10){$\dots$}}
\put(90,20){\circle{20}}
\put(156,53){\circle{20}}
\put(26,57){\circle{20}}
\put(160,20){\circle{20}}
\put(97,27){\line(2,1){49}}
\put(82,26){\line(-2,1){48}}
\put(24,46){\makebox(1,2){$\mag$}}
\put(155,42){\makebox(1,2){$\mag$}}
\put(159,9){\makebox(1,2){$\mag$}}
\put(100,20){\line(1,0){50}}
\put(26,57){\circle*{2}}
\put(160,20){\circle*{2}}
\put(156,53){\circle*{2}}
\put(150,20){\circle*{3}}
\put(90,20){\circle*{2}}
\put(170,20){\line(1,0){31}}
\put(210,20){\makebox(0,1){$\dots$}}
\put(217,20){\line(1,0){32}}
\put(337,2){\makebox(1,2){$(a_{1,1},b_{s,1})$}}
\put(260,0){\makebox(1,2){$(a_{1,1},\beta_s)$}}
\put(289,53){\makebox(1,2){$(a_{1,1},b_{s,2})$}}
\put(229,57){\makebox(1,2){$(a_{1,1},b_{s,i})$}}
\put(326,66){\makebox(6,10){$\mu_{s,2}$}}
\put(334,31){\makebox(6,10){$\mu_{s,1}$}}
\put(196,70){\makebox(6,10){$\mu_{s,i}$}}
\put(260,68){\makebox(6,10){$\dots$}}
\put(260,20){\circle{20}}
\put(326,53){\circle{20}}
\put(196,57){\circle{20}}
\put(330,20){\circle{20}}
\put(267,27){\line(2,1){49}}
\put(252,26){\line(-2,1){48}}
\put(194,46){\makebox(1,2){$\mag$}}
\put(325,42){\makebox(1,2){$\mag$}}
\put(329,9){\makebox(1,2){$\mag$}}
\put(270,20){\line(1,0){50}}
\put(196,57){\circle*{2}}
\put(330,20){\circle*{2}}
\put(326,53){\circle*{2}}
\put(320,20){\circle*{3}}
\put(260,20){\circle*{2}}
\end{picture}
\end{center}
We recall now the following definition and theorem from \cite{oka}.\\
\begin{DEF*}\cite{oka}
$$G_{m,n}:=<g_1, \ldots, g_m |\
g_k^{-1}(g_1\cdots g_n)g_{k+n}(g_1\cdots g_n)^{-1},\ \forall k=1\vir m>,$$
where the indices in the relators are taken to be cyclical mod $m$.
\end{DEF*}
\begin{teo}[\cite{oka}]\label{pi1dioka}
$$\pi_1(\ci^2\meno \{x^m=y^n\})\cong G_{m,n}$$.
\end{teo}
\begin{PROP}
If $B=\{f(x)=g(y)\}$ as above,
$$\pi_1(\ci^2\meno B)\cong G_{m,n}$$
where $n\ugu(n_1\vir n_r)$, $m\ugu(m_1\vir m_s)$ (the greatest common divisors).
\end{PROP}
\begin{DIM}
Every path in
$\Pi$ induces a braid (acting on $p^{-1}(a_{1,1})$) that is its
{\it braid monodromy}: we compute the relations in $\pi_1(\ci^2 \meno
B)$ by the braid monodromy of the generators of $\Pi$, following the
method introduced in \cite{moi}.\\
In order to express the braid monodromy of a path we use the standard 
generators of the braid group on $d$ strands
given by the positive half-twists $\sigma_i$, $1 \leq i \leq
d-1$ exchanging the $i$th and the $(i+1)$th strands
counterclockwise (the reader unexperienced with the braid group can
find precise definitions and more, e.g., in \cite{bir}).\\ 
In order to compute the braid monodromy of the generators we chose
for $\Pi$, we think the points $a_{i,h}$ lying on a line following the
lexicographical order in their indices, i.e. $a_{i,h}\mag a_{i',h'}$ $\iff$ 
$i\mag i'$ or, if $i\ugu i'$, $h\mag h'$.\\
The braid monodromy of $\rho_1$ is
$$
\tilde\sigma_1^{n_1}\cdots\tilde\sigma_s^{n_1}
$$
where $\tilde\sigma_j\ugu\sigma_{m_0+\cdots+m_j-1}\cdots\sigma_{m_0+\cdots+m_{j-1}+1}$,
($m_0\ugu0$), and gives us the relations
$$
\mu_{j,k}=T_{j;1,n_1}\mu_{j,k+n_1}T_{j;1,n_1}^{-1}
$$
for all $j,k$, where $T_{j;1,l}=\mu_{j,1}\mu_{j,2}\cdots\mu_{j,l}$ and the second
index of the $\mu_{j,k}$'s is taken to be cyclical (mod$m_j$).\\
Since by condition \ref{arba}, lifting the path $l_i$ gives the identity braid for all $i$,
then the braid monodromy of $\rho_i$ is similar, i.e.
$$
\tilde\sigma_1^{n_i}\cdots\tilde\sigma_s^{n_i}
$$
inducing in $\pi_1(\ci^2\meno B)$ (and  $\pi_1(\PP^2\meno\bar B)$),
the relations
$$
\mu_{j,k}=T_{j;1,n_i}\mu_{j,k+n_i}T_{j;1,n_i}^{-1}
$$
for all $i,j,k$. It is easy to see (cf. proposition 1.1 in \cite{mapi}) that these relations
are equivalent to
$$
\mu_{j,k}=T_{j;1,n}\mu_{j,k+n}T_{j;1,n}^{-1}
$$
for all $j,k$, where $n=(\unoenne nr)$.\\
The monodromy of $\Delta_{j;1,1}$ is retrieved from the braid $f^{-1}(\omega_j)$ and gives
the positive half twist
\begin{center}
\begin{picture}(300,120)(-20,40)
\put(-20,98){$\cdots$}
\put(10,100){\circle*{3}}
\put(5,92){$a_{j,1}$}
\put(40,100){\circle*{3}}
\put(35,92){$a_{j,2}$}
\put(90,100){\circle*{3}}
\put(85,92){$a_{j,m_j}$}
\put(130,100){\circle*{3}}
\put(116,112){$a_{j+1,1}$}
\put(60,98){$\cdots$}
\put(145,97){$\cdots$}
\put(170,100){\circle*{3}}
\put(156,112){$a_{j+1,\left[\frac{m_{j+1}}2\right]}$ }
\put(200,100){\circle*{3}}
\put(186,70){$a_{j+1,\left[\frac{m_{j+1}}2\right]+1}$}
\put(220,97){$\cdots$}
\put(260,100){\circle*{3}}
\put(246,112){$a_{j+1,m_{j+1}}$}
\put(280,97){$\cdots$}
\put(208,80){\vector(-1,3){6}}
\qbezier(10,100)(60,140)(110,100)
\qbezier(110,100)(145,65)(200,100)
\end{picture}
\end{center}
i.e.
\begin{equaz}
T^{-1}\sigma_{m_0+\cdots+m_j}T\label{bmo}
\end{equaz}
$$
T=(\sigma_{m_0+\cdots+m_j+1}\cdots\sigma_{m_0+\cdots+m_{\left[\frac{j+1}2\right]}})
(\tilde\sigma_j)
$$
and gives the relation
$$
\mu_{j,1}=(\mu_{j+1,1}\cdots\mu_{j+1,\left[\frac{m_{j+1}}2\right]})
\mu_{j+1,\left[\frac{m_{j+1}}2\right]+1}
(\mu_{j+1,1}\cdots\mu_{j+1,\left[\frac{m_{j+1}}2\right]})^{-1}.
$$
This relation is best understood in terms of the ``minimal standard
generators'' (cf. \cite{mapi})
$$
\gamma_{j,k}=(\mu_{j,1}\cdots\mu_{j,k-1})\mu_{j,k}^{-1}(\mu_{j,1}\cdots\mu_{j,k-1})^{-1}.
$$
The relation above becomes the simpler
$$
\gamma_{j,1}=\gamma_{j+1,\left[\frac{m_{j+1}}2\right]+1}.
$$
The braid monodromies of the other $\Delta_{j;i,h}$ are a conjugate of \ref{bmo} by a multiple
of $\tilde\sigma_j\tilde\sigma_{j+1}$ and give the relations
$$
\gamma_{j,k}=\gamma_{j+1,\left[\frac{m_{j+1}}2\right]+k}
$$
for all $j,k$.\\
These are cancellation relations, since we can express each $\mu_{j,k}$ in terms of the
$\mu_{1,k}$'s, and moreover give us the relations
$$
\mu_{1,k}=\mu_{1,k+m_j}
$$
for all $j,k$.\\
So, if $m:\ugu(m_1\vir m_s)$,  the paths $\mu_1\ugu\mu_{1,1},$ $\cdots$ $\mu_m\ugu\mu_{1,m},$ 
generate $\pi_1(\ci^2\meno B)$, and between them we have only
the relations  
$$
\mu_k=T_{1,n}\mu_{k+n}T_{1,n}^{-1}
$$
where $T_{1,n}=\cunoenne{\mu}n$ with cyclical indices mod $m$.
\end{DIM}
\begin{rem}\label{local-global}
By theorem \ref{pi1dioka}, $\pi_1(U_{i,j}\meno B)\cong G_{n_i,m_j}$
and $\unoenne{\mu}{m_j}$ are
(conjugated to) the standard generators for this group: in
particular, the map $\pi_1(U_{i,j}\meno B)\rightarrow\pi_1(\ci^2\meno B)$ induced by the inclusion
coincides with the map $(f_{\frac{n_i}n,\frac{m_j}m})_*$ we introduced
immediately after remark \ref{R}.\\
This implies that if $B$ is the branch curve of a normal generic cover
with monodromy $\mu:G_{m,n}\rightarrow\sgo_d$, the graph representing
the local monodromy at $P_{i,j}$ 
is the pullback by $(f_{\frac{n_i}n,\frac{m_j}m})$ of the graph
representing the global mo\-no\-dro\-my $\mu$.
\end{rem}
\begin{PROP}\label{pi1complproj}
$$\pi_1(\PP^2\meno\bar B)\cong\modulo{G_{m,n}}{<({\mu_1 \cdots \mu_m})^{\frac dm}>}$$
where $n\ugu(n_1\vir n_r)$, $m\ugu(m_1\vir m_s)$.\hfill\CVD
\end{PROP}
\begin{DIM}
To compute $\pi_1(\PP^2\meno B)$ we use the standard remark that the
kernel of the surjective map 
$$\pi_1(\ci^2\meno B)\rightarrow\pi_1(\PP^2\meno\bar B)\rightarrow0$$
is infinite cyclic and is generated by a loop $L$ around the line at infinity.\\
In our case this loop is
$$
L=(\mu_{1,1}\cdots\mu_{1,m_1})(\mu_{2,1}\cdots\mu_{2,m_2})\cdots(\mu_{s,1}\cdots\mu_{s,m_s})
$$
or in terms of the generators $\mu_i$
$$
L=(\mu_1\cdots\mu_m)^{\frac{m_1}m}\cdots(\mu_1\cdots\mu_m)^{\frac{m_s}m}=
(\mu_1\cdots\mu_m)^{\frac dm}
$$
\end{DIM}
Assume now $\bar B$ irreducible, i.e. $(n,m)\ugu1$.\\
In this case the monodromy of
the cover lifts to a {\rm generic} (geometric loops map to
transpositions) homomorphism 
$\mu:G_{n,m}\rightarrow\sgo_d$ for which $\mu(T_{1,n}^{\frac
  dn})=1$. By the classification of generic homomorphisms in
theorem \ref{classgraf}, the
mo\-no\-dro\-my graph is (exchanging $n$ and $m$ if necessary) a polygon.\\
We know that in this case there exist $h,k,a,b$ s.t. $n\ugu a(h\piu k)$ and $m\ugu bkh$ with $(h,k)\ugu 1$.
Now we can introduce our family: we define
$$
\bar g_l(x,w)=(x-w)(x-2w)\cdots(x-lw)
$$
$$
\bar f_l(y,w)=(y-w)(y-2w)\cdots(y-lw)
$$
and, given $h,k$ coprime, we consider the generic cover of degree
$h+k$ branched over
$$
\bar g_{h+k}(x,w)^{hk}=\bar f_{hk}(y,w)^{h+k}
$$
with monodromy graph a polygon with $m$ edges, valence $1$ and
increment $h$.\\
Here all the singularities have the same form $x^{hk}=y^{h+k}$, so, by remark 
\ref{local-global}, all the local monodromy graphs have to coincide
with the global one.\\
In order to ensure the existence of the cover we have only to check
that the monodromy of $(\mu_1 \cdots \mu_{h+k})^{hk}$ is trivial, which
was clear since the very beginning because it belongs to the center of
the (local) fundamental group (in fact, the order of the monodromy of $\mu_1 \cdots
\mu_{h+k}$ is exactly $hk$).\\
Moreover, by corollary \ref{smoothness}, having all the singular
points $a=b=1$ the surface we defined is smooth, whence the cover is
smooth if and only if $h=1$; in this last case one can easily check
%
that the cover is given by the projection on the plane $z=0$ from the
point $(0,0,1,0)$ of the surface
$$
z^{k+1}-(k+1)z\bar f_{k}(x,w)+k\bar g_{k+1}(y,w)=0
$$
Finally we can state the counterexample we were looking for:
\begin{PROP}
Let $B$ be the projective plane curve of degree $30$ given by the equation
$$
\bar g_{5}(x,w)^{6}=\bar f_{6}(y,w)^{5},
$$
Then there are two generic covers $S' \stackrel{\pi'}{\rightarrow}
\PP^2$, $S'' \stackrel{\pi''}{\rightarrow} \PP^2$, with 
\begin{enumerate}
\item\label{sonolisci} $S'$, $S''$ smooth;
\item\label{grado6} deg $\pi'=6$;
\item\label{grado5} deg $\pi''=5$;
\item\label{unoliscio} the ramification divisor of $\pi'$ is smooth;
\item\label{2singolare} the ramification divisor of $\pi''$ has
  (exactly) $30$ ordinary cusps as singularities.
\end{enumerate}
\end{PROP}
\begin{DIM}
The covers $\pi'$, $\pi''$ are the covers of the family we just
constructed for, respectively, $h=1, k=5$ and $h=2, k=3$. 
We check quickly the $5$ properties:\\
(\ref{sonolisci}) holds for every surface in our family;\\
(\ref{grado6}) and (\ref{grado5}) follows because the degree is the number of
vertices of the graph, i.e. $h+k$.\\
Finally, (\ref{unoliscio}) and (\ref{2singolare}) come directly from remark
\ref{R}.
\end{DIM}
This is a counterexample to the Chisini's conjecture if we drop the
assumption that the ramification divisor is non-singular.\\
This family does not produce counterexamples in
higher degrees: in fact the pair $(5,6)$ is the only one which can
be expressed as sum and product of two coprime integers in two different ways.\\
Indeed, suppose we have $h\piu k\ugu h'k'$ and $hk\ugu h'\piu k'$ with $(h,k)=(h',k')=1$ and, say,
$h\min k$, $h'\min k'$, $h\piu k\min hk$.\\
From $h'\piu k'\mag h'k'$ we have that $h'\ugu1$ and $k'\ugu hk\men1$.\\
But now $k(h\men1)\ugu h\piu1\ugu h\men1\piu2$ and it must be $(h\men1)|2$,
so that $h\ugu 2$ and $k\ugu 3$
which gives $k'\ugu
5$.\\
In order to find counterexamples to a Chisini-Kulikov-Nemirovski's type
result in arbitrarily large degrees, we have to consider a
slightly different family:
\begin{PROP}
Let $t \in \eNNe$, $t\geq 2$ 
$B$ be the projective plane curve given by the equation
$$
\bar g_{4t+1}(x,w)^{2t(2t+1)}=\bar f_{2t(2t+1)}(y,w)^{4t+1}.
$$
Then there are two generic covers $S' \stackrel{\pi'}{\rightarrow}
\PP^2$, $S'' \stackrel{\pi''}{\rightarrow} \PP^2$, with  
$S'$, $S''$ smooth, degrees respectively $4t+1$ and $4t+2$, 
and both singular ramification divisor.
\end{PROP}
In fact, the case $t=1$ is exactly the case of previous proposition,
so the statement still holds except for the singularities of the ramification divisor.\\
\begin{DIM}
The cover of degree $4t+1$ is the cover in our family for $h=2t$,
$k=2t+1$.\\
The cover of degree $4t+2$ is simply the cover constructed in the same
way as we did for our family, starting from the monodromy graph given
by the polygon with $4t+2$ vertices, valence $t$ and increment $1$.\\
The smoothness comes from corollary \ref{smoothness} observing that
locally we have $h=b=1$. \\
The other verifications are exactly as in the
previous case and we leave them to the reader.
\end{DIM}
%
%

%
%
%
\end{document}